\documentclass[12pt]{amsart}
\usepackage{amssymb,amsmath,amsthm,amscd}
\usepackage[all]{xy}

\usepackage[hypertex]{hyperref}

\numberwithin{equation}{section}

\newtheoremstyle{fancy1}{10pt}{10pt}{\itshape}{12pt}{\textsc\bgroup}{.\egroup}{8pt}{
}
\newtheoremstyle{fancy2}{10pt}{10pt}{}{12pt}{\itshape}{.}{8pt}{ }

\theoremstyle{fancy1}

\newtheorem{cor}[equation]{Corollary}

\newtheorem{prop}[equation]{Proposition}
\newtheorem{thm}[equation]{Theorem}

\newtheorem*{main*}{Theorem}

\newtheorem*{cor*}{Corollary}

\newtheorem*{problem*}{Problem}

\theoremstyle{fancy2}

\newtheorem*{rem*}{Remark}

\newcommand{\cref}[1]{Corollary~\ref{#1}}

\newcommand{\CP}{\mathbb{C\mkern1mu P}}

\newcommand{\Sph}{\mathbb{S}}

\newcommand{\R}{{\mathbb{R}}}

\newcommand{\MV}{{\mathcal{V}} }
\newcommand{\MH}{{\mathcal{H}} }

\newcommand{\SO}{\ensuremath{\operatorname{SO}}}

\newcommand{\SU}{\ensuremath{\operatorname{SU}}}

\newcommand{\T}{\ensuremath{\operatorname{T}}}
\renewcommand{\S}{\ensuremath{\operatorname{S}}}

\newcommand{\fg}{{\mathfrak{g}}}
\newcommand{\fk}{{\mathfrak{k}}}

\newcommand{\fm}{{\mathfrak{m}}}

\def\con#1=#2(#3){#1 \equiv #2 \bmod{#3}}

\newcommand{\To}{\longrightarrow}

\renewcommand{\sec}{\ensuremath{\operatorname{sec}}}

\DeclareMathOperator{\Id}{Id}

\newcommand{\no}{\noindent}

\newcommand{\spa}{\mbox{span}}

\newcommand{\nnc}{non-negative curvature}

\newcommand{\norm}[1]{\left\Vert#1\right\Vert}
\newcommand{\abs}[1]{\left\vert#1\right\vert}

\def\bsm{\begin{smallmatrix}}
\def\esm{\end{smallmatrix}}
\def\bpm{\begin{pmatrix}}
\def\epm{\end{pmatrix}}
\def\beq{\begin{equation}}
\def\eeq{\end{equation}}

\begin{document}

\begin{titlepage}
\center{\Large{{ On
M. Mueter's Ph.D. Thesis on Cheeger deformations}}}\\[5cm]
\center{\Large{{ Lectures by Wolfgang Ziller }}}\\[1cm]
\center{\Large{{ University of Pennsylvania,  Dec. 18, 2006 }}}\\[1cm]
\center{\Large{{ Notes by Chenxu He }}}\\[3cm]
\end{titlepage}

In 1973 Cheeger \cite{C}, motivated by the Berger metrics on spheres, studied
a general construction for a Riemannian manifold with an isometric
group action that shrinks the metric in the direction of the orbits.
If the original metric has \nnc, so does the new one and it tends to
have less 2-planes with 0 curvature. In 1987 Michael M\"uter studied this
construction and its curvature  properties in detail in his Ph.D.
thesis \cite{M} at the University of M\"unster under the direction of
Wolfgang Meyer. One can view it as a deformation of the original
metric (although this was not done so in \cite{C}), and we will call it, following M\"uters terminology, a {\it
Cheeger deformation}. His Ph.D. thesis, which was written in German, was never published or
translated. Although there is no main quotable theorems in his
thesis, it nevertheless contains a wealth of information about such
Cheeger deformations. Since they have been used extensively over the
last 20 years (see \cite{Z} for a survey of examples obtained by using Cheeger deformations), I  gave some talks in  my secret seminar at the
University of Pennsylvania in 2006 summarizing the results.  Since this may be of
some interest to others as well, I decided it might be useful to post it on the arxive.
 A scanned copy of M\"uter's thesis can be obtained on my home page at
 www.math.upenn.edu/\~{}wziller/research.

\bigskip

 Suppose $(M,g)$ is a Riemannian manifold and $G$  a
Lie group acting isometrically on $M$. Let
$\pi\colon M\to M/G$ be the projection, which is a Riemannian submersion
on the regular part of the group action. We assume that $G$ is
compact and choose a fixed bi-invariant metric $Q$ on the Lie
algebra $\fg$. On $M \times G$, we  define the metric $g +
\frac{1}{t}Q$, and consider a second (smooth) Riemannian  submersion
\begin{center}$
\begin{array}{rcccl}
\sigma & : & M \times G & \To & M \\
    &   & (p,g) & \longmapsto & g^{-1}p
\end{array}$
\end{center}
which induces a metric $g_t$ on $M$. It extends smoothly across
$t=0$ such that $g_0=g$ and we thus obtain a variation of $g$. One can also consider $g_t$ for $t<0$ since
for small negative $t$ the metric $g_t$ is still positive definite.
 If the original metric has \nnc,
so does $g_t$ for $t>0$ by  O'Neill's formula. This is usually not the case for $t<0$.

\smallskip

 For each fixed $p
\in M$, let $G_p$ be the isotropy group with Lie algebra $\fg_p$,
and define the $Q$-orthogonal decomposition $\fg = \fg_p + \fm_p$.
Using the action fields $X^*$ on $M$,\; $X\in\fg$ , we identify
 $\fm_p$, via $
 X\in\fm_p\to X^*(p)$, with the tangent space of the orbit $Gp$ at $p$.

 At $p \in M$, we  define the following two subspaces of $T_pM$:
\begin{eqnarray*}
\mathcal{V}_p & = & T_pG(p) = \{X^*_p \mbox{ } |\mbox{ } X \in \fg\} \\
\mathcal{H}_p & = & \{\mbox{ } \xi \in T_pM \mbox{ }| \mbox{ } g_t(\xi, \mathcal{V}_p) = 0\}
\end{eqnarray*}
We call $\mathcal{V}_p$ and $\mathcal{H}_p$ the \emph{vertical
space} and \emph{horizontal space} at $p$, although in general these
are not smooth distributions. For $ V\in T_pM$, we denote by
$V^{\mathcal{V}}$ and $V^{\mathcal{H}}$ the vertical and horizontal
component of $V$. For any $V \in T_pM$, there exists a unique $X \in
\fm_p$, such that $V^{\mathcal{V}} = X^*_p$,  and we denote this
vector $X$ by $V_\fm$, i.e., $V_\fm^*(p)=V^{\mathcal{V}}_p$.

\smallskip

 We now
define an \emph{orbit tensor} $P_t\colon\fm_p\to\fm_p$ at a point
$p$ as
\[
Q(P_t(X),Y) = g_t(X^*,Y^*)_p, \mbox{  }\forall \mbox{ } X, Y \in
\fm_p.
\]
For $t=0$ we set for simplicity $P_0=P$.

\no We define another tensor $C_t\colon T_pM\to\T_pM$ at $p$, called
the \emph{metric tensor} by
\[
 g(C_t(X),Y)=g_t(X,Y), \mbox{  } \forall \mbox{ } X, Y \in T_pM.
\]
\no Thus $C_t=\Id$ on $\MH_p$ and $C_t=P^{-1}P_t$ on $\MV_p\simeq  \fm_p$.

\begin{prop} The orbit tensor and metric tensor are given by
\begin{enumerate}
\item[(a)] $P_t  =  (P^{-1} + tI)^{-1} = P(I + tP)^{-1}$,
\item[(b)] $C_t(V)  =  (I + tP)^{-1}(V^{\mathcal{V}}) + V^{\mathcal{H}}$.
\end{enumerate}
where $P=P_0$.
\end{prop}

Indeed, under the Riemannian submersion $\sigma$, the horizontal lift of $X=
X^\mathcal{V} + X^\mathcal{H}=X^*_\fm + X^\mathcal{H} \in T_pM$ is given by
\begin{equation*}\label{lift}
(P^{-1}(P^{-1}+t\Id)^{-1}X_\fm^*+X^\mathcal{H},-t(P^{-1}+t\Id)^{-1}X_\fm)\in\T_pM\times\fg
,\end{equation*}
and thus
\begin{align*}
Q((P^{-1}+t\Id)^{-1}X_\fm^*,P^{-1}(P^{-1}+t\Id)^{-1}X_\fm^*)&+\frac 1 t
Q(t(P^{-1}+t\Id)^{-1}X_\fm,t(P^{-1}+t\Id)^{-1}X_\fm)\\
&=Q((P^{-1}+t\Id)^{-1}X_\fm,X_\fm)
\end{align*} on $\MV$, whereas
 orthogonal to $\MV$ the metric is unchanged.

Now we consider the sectional curvature in the metric $g_t$. M\"uter made the crucial observation that
 instead of describing the change of the sectional curvature of a fixed plane spanned by $X,Y$ in the metric $g_t$,  it is better to consider the change of the curvature of the moving  planes $\spa\{C_t^{-1}X,C_t^{-1}Y\}$. Indeed:
 \begin{equation}\label{lift} \text{ The horizontal lift of } C_t^{-1}X = (\Id+tP)X_\fm^*+X^\MH \text{ is equal to }
 (X,-tPX_\fm).\end{equation}
  This implies that if the plane $\spa\{X,Y\}$ has non-negative curvature in $g_0$, then so do the planes $\spa\{C_t^{-1}X,C_t^{-1}Y\}$ in $g_t$. This is not necessarily the case for the fixed plane $\spa\{X,Y\}$, in fact its curvature can become negative for some $t>0$. Thus 0-curvature planes tend to move with $C_t^{-1}$.  If we let
$$
\kappa_c(t)  = g_t(R_t(C_t^{-1}(V),C_t^{-1}(W))C_t^{-1}(V),C_t^{-1}(W))
$$
be the un-normalized
sectional curvature, we obtain:
\begin{prop} The sectional curvature of the metric $g_t$ is determined by
$$\kappa_c(t)=  g(R(V,W),V,W) + \frac{t^3}{4}\norm{[P(V_\fm),
P(W_\fm)]}_Q^2 +
            z(V,W,t)$$
            where
            $$z(V,W,t) = 3 \mbox{ }t \max_{Z \in \fg, \abs{Z} =
1}\frac{\{dw_Z(V,W) + \frac{t}{2}Q([P(V_\fm),P(W_\fm)],Z)\}^2}{t
\mbox{ } g(Z^*, Z^*)+1}$$ and $w_Z(U) = g(Z^*, U)$ $\forall \mbox{
}U \in T_pM$ is the one form dual to $Z\in\fg$. Furthermore, if $X,Y$ are orthonormal in $g_0$,
 $$||C_t^{-1}(V) \wedge C_t^{-1}(W)||^2_{g_t}   =  t^2 \norm{[P(V_\fm), P(W_\fm)]}_Q^2 +\\
t(\norm{P(V_\fm)}_Q^2 + \norm{P(W_\fm)}_Q^2) +1.$$
\end{prop}

\smallskip
The first two terms are clearly the curvature in $M\times G$, and $z(V,W,t)$ is
the O'Neill term of the submersion $\sigma$. Hence

\begin{cor} If $g_0$ has non-negative curvature and $V,W$ span a 0 curvature plane in $g_o$, then
\begin{enumerate}
\item[(a)]  $ \kappa'_c(0) \geq 0$, and $ \kappa'_c(0) = 0$ if and only if $dw_Z(V,W) = 0, \forall \mbox{ } Z \in \fg$.
\item[(b)] If $ \kappa'_c(0) = 0$ and $
[P(V_{\fm}), P(W_{\fm})] \ne 0$  then $ \kappa''_c(0) = 0$,\, $
\kappa'''_c(0) > 0$ and  $\kappa_c(t) > 0$ for all $t>0$.
\item[(c)] If $ \kappa'_c(0) = 0$ and $
[P(V_{\fm}), P(W_{\fm})] = 0$  then   $\kappa_c(t) = 0$ for all
$t>0$.
\end{enumerate}
\end{cor}
In particular, if $
[P(V_{\fm}), P(W_{\fm})] \ne 0$, i.e., the plane $\spa\{
P(V_{\fm}), P(W_{\fm})\}$ has positive curvature in the biinvariant metric $Q$ on $G$, then $\sec(C_t^{-1}X,C_t^{-1}Y)>0$ for all $t>0$. Furthermore, either $\sec(C_t^{-1}X,C_t^{-1}Y)>0$ for all $t>0$, or $\sec(C_t^{-1}X,C_t^{-1}Y)=0$ for all $t>0$.

\smallskip

M\"uter observes that in a product metric, where mixed 2-planes have 0 curvature, $k_c'(0)=0$ for such 2-planes and hence they become positive for a Cheeger deformation (which we can assume  acts  on both factors) iff $
[P(V_{\fm}), P(W_{\fm})] \ne 0$.

\smallskip

One  obtains an obstructions to deforming non-negative curvature to positive curvature with a Cheeger deformation:

\begin{cor} If $g$ has non-negative curvature, a 2-plane tangent to a totally geodesic flat 2-torus which contains horizontal vectors
remains flat.
\end{cor}
Indeed, by applying the Gauss Bonnet formula to the metrics $g_t$ restricted to the 2-torus, one sees that $ \kappa'_c(0) = 0$, and since the 2-plane contains horizontal vectors $
[P(V_{\fm}), P(W_{\fm})] = 0$.

\smallskip

For the values of $dw_Z$, which determine $k_c'(0)$, we have:
\begin{prop}Given $Z\in\fg$,
$dw_Z$ is determined by:
\begin{enumerate}
\item[(a)]  $2dw_{Z}(X^{*},
Y^{*}) = Q([P(X),Y] + [X,P(Y)] - P([X,Y]), Z) $ for $X,Y\in \fm\subset \fg$.
\item[(b)] $2dw_{Z}(X^{*},
U)=-g(\nabla_{X^*}U,Z^*)=-2g(S_U(X^*),Z^*) =-U(Q(PX,Z)$ for  $X\in \fm$ and $U\in \MH$.
\item[(c)] $dw_{Z}(U,V) =-g(\nabla_UV,Z^*)= -g(A_UV,Z^*)$ for $U,V\in \MH$
\end{enumerate}
\end{prop}
Here $S_U$ is the shape operator of the orbit in direction of the normal $U$, and $A_UV$ the O'Neill tensor of the submersion $\pi$, at least on the regular part. To extend these tensors, and thus Proposition 1.5,  to smooth tensors on all of $M$ one defines them as follows:

$$
g(S_U(X^*),Y^*)=g(\nabla_{X^*}Y^*,U) \text{ and } g(A_UV,X^*)=-dw_X(U,V) \text{ with } X,Y\in\fm, U,V\in\MH
$$
Notice that $A$ cannot be defined on all of $M$ as usual since $U\in\MH_p$ may not be extendable to a smooth horizonal vector field near $p$.

\begin{cor} Assume that $g_0$ has non-negative curvature. Then:
\begin{enumerate}
\item[(a)]  a horizontal 2-plane $\spa\{C_t^{-1}U,C_t^{-1}V\}=\spa\{U,V\}$ has 0-curvature in the metric $g_t$ iff $\spa\{U,V\}$ has 0-curvature in $g$ and $A_UV=0$.
\item[(b)]   a vertizontal 2-plane  $\spa\{C_t^{-1}X^*,C_t^{-1}U\}=\spa\{C_t^{-1}X^*,U\}$ has 0-curvature in the metric $g_t$ iff $\spa\{X^*,U\}$ has 0-curvature in $g$ and $S_U(X^*)=0$.
\end{enumerate}
\end{cor}

\smallskip

M\"uter observes that the $G$ action on the right of the second factor in $M\times G$ descends under $\sigma$, to an isometric action in $g_t$. In general, $C(G)\times G$ still acts by isometries in $g_t$, where $C(G)$ is the centralizer of $G$ in the full isometry group of $g_0$. He also examines properties of Cheeger deformed metrics on quotients $M/L$ where $L\subset C(G)$ or $L\subset G$ acts freely on $M$.

\smallskip

In particular, if $ K_1 \subset K_2 \subset \cdots \subset K_n \subset G$ is a chain of subgroups, one can iterate a Cheeger deformation one step at a time by a paremeter $t_i$. He shows that such a metric takes on the form
$
\sum_i s_i Q_{\fk_i\cap \fk_{i-1}^\perp}
$
where $s_i=(1+t_1+\dots + t_i)^{-1}$. He then derives general criteria under which a 2-plane in this metric has 0 curvature.

As a special case, he discusses in more detail the case $n=1$ and shows that the curvature of
 $g_s= sQ_{|\fk}+Q_{|\fm}$, $\fm=\fk^\perp$, is given by:

\begin{align}
g_s(R_s(A, B)A,B) &=  \tfrac{1}{4}\left\|s[A,B]_\fm + (1-s)[A_\fm,B_\fm]_\fm \right\|^2_Q +
    (1-s)^2\left\|[A_\fk,B_\fk] \right\|^2_Q  +\\
  &-(1-s)(2-s)Q([A_\fk,B_\fk] , [A,B]_\fk) -
 \tfrac{4-3s}{4}\left\|[A,B]_\fk \right\|^2_Q\notag
\end{align}

He observes that this implies that in general such metrics only have non-negative curvature when $s\le 1$.

\smallskip

In the last part, M\"uter considers two sided  iterated Cheeger deformations on a compact Lie group $G$
with biinvariant metric $Q$. Let $H$ be a closed subgroup of $G$ which acts on $G$ by
right translations. We consider two sequences of closed subgroups of $G$:
\[
 K'_1 \subset K'_2 \subset \cdots \subset K'_n \subset G \supset K_n \supset \cdots
\supset K_2 \supset K_1 \supset H
\]
$K'_i$ acts on $G$ from the left and $K_i$ acts from the right, $i = 1 \ldots n$. M\"uter considers three types of metrics:\\

\textit{Type 1.} \emph{The initial metric on $G$ is left-invariant and $H$ right-invariant.
The new metric on $G/H$ is obtained by applying an iterated Cheeger deformation by $K'_i$ from the left.}\\

\textit{Type 2.} \emph{The initial metric on $G$ is obtained from a biinvariant metric on $G$ by iterated Cheeger deformations with $K_i$ from the right. The new metric on $G/H$ is obtained by applying to this metric the iterated
Cheeger deformation by $K'_i$ from the left.}\\

\textit{Type 3.} \emph{A metric on $G$ is obtained from a right-invariant metric by
applying an iterated Cheeger deformation by $K_i$ from the right.  The  metric on $G/H$ is the submersed metric. }\\

Notice that these metrics on $G/H$ are typically not homogeneous anymore.
In order to study these metrics, M\"uter proves the following theorem on biquotients of product Lie groups:
\begin{thm}
Suppose $G_1$, $G_2$ are two compact Lie groups, $g_1$ is a left-invariant metric on $G_1$ and
$g_2$ is a biinvariant semi-Riemannian metric on $G_2$. Let $G = G_1 \times G_2$ and $g = g_1 +
g_2$ a product metric on $G$. Assume that $U \subset G \times G$ acts freely on $G$ with $U_R \subset
G_2$ where $U_R$ is the projection of $U$ onto the second component. If the biquotient $(G/\!/U, g)$ is a Riemannian manifold with
positive sectional curvature, then any biinvariant metric on $G$ induces a metric   with
positive sectional curvature on $G/\!/U$.
\end{thm}

The following are all consequences of this general result, although not all are obvious.

\begin{cor}
An iterated Cheeger deformation of a symmetric metric on a symmetric space of rank $> 1$ (a \textit{Type 1 metric}) has
$0$-curvature planes.
\end{cor}
\begin{cor}
If $G/H$ with a \textit{Type 2} or \textit{Type 3} metric has positive sectional curvature, then
$G/H$ has a normal homogeneous metric with positive sectional curvature.
\end{cor}
\begin{cor}
 If a \textit{Type 1} metric on $G$ (i.e. $H=\{e\}$) has positive sectional curvature, then $G = SU(2)$ or $SO(3)$.
\end{cor}

M\"uter calls a homogeneous space $G/H$ half-irreducible if there exists a chain of subgroups $K_0=H\subset K_1\subset \dots \subset K_n=G$ such that the adjoint representation of $H$ induces on $\fk_{i+1}/\fk_i$ an irreducible representation for all $i$.

\begin{cor}
 If a half-irreducible homogeneous space admits a homogeneous metric $g$ such that some iterated Cheeger deformation of $g$ from the left has positive sectional curvature, then $G/H$ admits homogeneous metric with positive curvature.
\end{cor}

Some further discussions in his thesis:

\smallskip

1) A Killing vector field $X$ defines an isometric action by $\R$ on M and M\"uter considers a Cheeger deformation $g_t$ with respect to this action. He shows that if $X$ has unit length, then $g_t$ has curvature bounded independent of $0<t<1$, but has volume going to $0$ (vanishing minimal volume).
He then observes that under the condition that $X$ has no zeroes, one easily changes the metric so that $X$ is Killing with unit length, and thus $M$ has vanishing minimal volume as well.

\smallskip

2) He considers a Cheeger deformation of the product of two round sphere metrics in $\Sph^2\times\Sph^2$ by $\SO(3)$ acting diagonally. He shows that at every point not lying on a diagonal or anti-diagonal, there exists a unique 0 curvature plane, and at these exceptional points, a circle's worth of 0 curvature planes. He also remarks that, in contrast, under a Cheeger deformation of the non-negatively curved Cheeger metric on $\CP^2 \# -\CP^2$,
which can be regarded as $\Sph^3\times_{\S^1}\Sph^2$, by the isometric action of $\SU(2)$, all 0 curvature planes remain 0 curvature planes.

\smallskip

3) Furthermore, he gives a simple proof, using Cheeger deformations, of a theorem by Wallach \cite{W}  that produces all positively curved homogeneous metrics.

\bigskip

This is not intended to be a complete list of all his results, but a selection of the ones that seemed most interesting to me. All misrepresentations  of M\"uter's results (and possible  mistakes) are to be blamed on this author of course.

\bigskip

\begin{center}
\bf More recent Applications
\end{center}

Cheeger deformations have been used frequently over the last 20 years, and thus some results
 in M\"uter's thesis have been rediscovered. Usually without reference though, since the thesis
  was only available to a few. We collect here a some of these applications.

In addition to the homogeneous examples, all known biquotient metrics with positive curvature
\cite{E1,E2,E3,B} are obtained as Cheeger deformations of biinvariant metrics (see  \cite{Z} for a survey).
In\cite{D} curvature properties of iterated Cheeger deformed metrics on Eschenburg spaces were examined.

Most of the known examples of metrics which are almost or quasi positively curved are also obtained
via a Cheeger deformation (sometimes iterated) of a biinvariant metric, see \cite{PW1,W,Wi,EK,Ke1,Ke2}.

In \cite{GZ} it was shown that  a "scaled up" Cheeger deformation (i.e. $t<0$) of a biinvariant
metric in the direction of an abelian subalgebra has nonnegative curvature for $t>-1/4$. Although
this is an immediate consequence of his formula (1.7), this fact escaped his attention. It was used
in \cite{GZ} to  construct many non-negatively curved metrics on cohomogeneity one manifolds, giving
rise to non-negative curvature on 7 dimensional exotic spheres.  In \cite{Sch} it was shown that
on a simple Lie algebra only abelian subalgebras can be scaled up and still remain non-negatively
curved. On the other hand, in \cite{ST} the case of homogeneous fibrations $K/H\to G/H\to G/K$ was
studied and it was shown that a scaling up for the Cheeger deformation of $K$ on $G/H$
(i.e. enlarging the fibers $K/H$) has non-negative curvature if $(K,H)$ is a symmetric pair.

In \cite{HT} and \cite{ST} it was observed that the space of homogeneous metrics on $G/H$ with nonnegative
 curvature is connected, since  in a Cheeger deformation by the left action of $G$ on $G/H$
 the rescaled metric
$ t g_t$ has non-negative curvature  and clearly converges to the biinvariant metric $Q$ as $t\to \infty$.
K.Tapp calls this an inverse linear path (one adds $t\Id$ to the inverse) and examines the Taylor
series of this path at the biinvariant metric $Q$ to examine which initial directions give rise to a
linear path with non-negative curvature. So far, among left invariant metrics on compact Lie groups,
the only known examples with non-negative curvature are obtained by  combing iterated Cheeger
deformations with scaling up in the direction of abelian subalgebras.

 \eqref{lift} immediately implies the following fact: If one uses a Cheeger deformation with a
 Lie group $G=\SU(2)$ or $\SO(3)$, then any 2-plane that has a two dimensional projection onto a
  $G$ orbit becomes positively curved for $t$ large, even if the plane initially has negative
  curvature.  This was used in \cite{PW2} in some of their  deformations in their
  recently proposed construction of a metric of
  positive curvature on the Gromoll-Meyer sphere \cite{GM}.

As B.Wilking observed, one obtains a simple proof of the Schwachh\"ofer-Tuschmann theorem \cite{ST}
that any cohomogeneity one $G$-manifold carries an almost non-negatively curved metric. One puts a
metric with non-negative curvature near the singular orbits (since these are homogeneous disc bundles)
 and extends arbitrarily in the middle. In a Cheeger deformation with $G$, all curvatures in the middle
  go to 0 as $t\to\infty$. Indeed, this follows immediately from Proposition 1.2 since any 2-plane has a
   non-trivial projection onto a $G$ orbit.

As F.Wilhelm and W.Ziller observed, one also obtains a simple proof of the Fukaya-Yamaguchi \cite{FY}
theorem that
the total space of a principal $G$-bundle (with $G$ compact)  over an
almost non-negatively curved manifold is  almost non-negatively
curved as well. Indeed, one can put a connection metric on the total space which induces the given
almost non-negatively curved metric on the base. As above, Cheeger deformation by $G$ with $t\to\infty$,
produces the desired metric. For a 2-plane with non-trivial projection onto the fibers, this follows as
before. For a horizontal 2-plane the O'Neill term  goes to zero since the metric on the fiber goes to 0,
and hence the curvature converges to that of the base by O'Neill's formula.  The original proofs in
\cite{ST} and \cite{FY} involved quite complicated curvature calculations.

\providecommand{\bysame}{\leavevmode\hbox
to3em{\hrulefill}\thinspace}

\end{document}